\documentclass[3p,preprint,times]{elsarticle_preprint}

\usepackage{amsmath}
\usepackage{amsfonts}
\usepackage{amssymb,latexsym}
\usepackage{amsthm}
\usepackage{graphicx}

\newtheorem{remark}{Remark}[section]

\usepackage{color}
\usepackage{newlfont}
\usepackage{subfigure}
\usepackage{verbatim}
\usepackage{rotating}
\usepackage{multirow}
\usepackage{units}
\usepackage{bm}
\usepackage[english]{babel}
\usepackage[utf8]{inputenc}
\usepackage{algorithm}
\usepackage[noend]{algpseudocode}
\usepackage{booktabs}
\usepackage{caption}
\usepackage{hyperref}

\newcommand{\RR}{\mathbb{R}}

\journal{R. Cavoretto, A. De Rossi}

\begin{document}

\begin{frontmatter}

\title{Adaptive meshless refinement schemes for RBF-PUM collocation}


\author[address-TO]{R. Cavoretto} 
\ead{roberto.cavoretto@unito.it}

\author[address-TO]{A. De Rossi}
\ead{alessandra.derossi@unito.it}


\address[address-TO]{Department of Mathematics \lq\lq Giuseppe Peano\rq\rq, University of Torino, via Carlo Alberto 10, 10123 Torino, Italy}

\begin{abstract}
In this paper we present an adaptive discretization technique for solving elliptic partial differential equations via a collocation radial basis function partition of unity method. In particular, we propose a new adaptive scheme based on the construction of an error indicator and a refinement algorithm, which used together turn out to be ad-hoc strategies within this framework. The performance of the adaptive meshless refinement scheme is assessed by numerical tests.
\end{abstract}




\end{frontmatter}


\section{Introduction}\label{sec1}

Meshless methods are known to be powerful computational tools for solving approximation problems, which include either multivariate data interpolation or numerical solution of partial differential equations (PDEs). In literature a very popular approach to face this type of problems is based on the use of radial basis functions (RBFs). In fact, the meshless nature of RBF methods guarantees a sure flexibility with respect to geometric problem, a certain simplicity in the implementation in higher dimensions, and a high convergence order of the numerical method (see \cite{fas07,wen05}). 

In this paper we are interested in constructing adaptive discretization techniques for elliptic PDEs such as a Poisson problem. We test our adaptivity on the RBF partition of unity method (RBF-PUM). The basic idea of PUM consists of decomposing the original domain into several subdomains or patches forming a covering of it, then constructing a local RBF approximant on each subdomain. Such an approach, as evident in our numerical analysis, enables us to significantly reduce the computational instability due to (usual) severe ill-conditioning generated by RBFs. Originally, the PUM was introduced in \cite{bab97,mel96} to solve PDEs. Recently, this method has reached popularity because it allows one to efficiently and accurately solve large scale interpolation problems \cite{cav12,cav16,wen02}, as well as differential ones \cite{her16,lar17,saf15}. However, the problem of constructing a (really effective) adaptive method is still open and does not seem to have received the proper attention in literature. In doing so, our main focus is primarily addressed to find any \textsl{error indicator} and \textsl{refinement algorithm}, which can be used within a RBF-PUM collocation scheme. We remark that adaptive techniques have been intensively investigated for meshless RBF methods, but most of them mainly refer to either weak form and finite difference methods or collocation multiscale methods (see e.g. \cite{dav11,far13,oan17}). The study of new numerical tools for solving PDEs is further motivated by the fact that they govern many physical phenomena and models of applied mathematics and engineering, including for instance the distribution of heat/temperature, the propagation of sound or light waves, and fluid dynamics \cite{fas15}.    

The paper is organized as follows. In Section \ref{RBF-PUM} we focus on the RBF-PUM collocation scheme to solve a Poisson problem. In Section \ref{amrs} we explain the different phases of our adaptive refinement algorithm, also proposing some error indicators and a refinement strategy. In Section \ref{num_exp} we report numerical experiments devoted to illustrate the performance of our adaptive scheme. 


\section{RBF-PUM collocation for Poisson problems} \label{RBF-PUM}


Given an open and bounded domain $\Omega \subseteq \RR^s$ partitioned into $d$ subdomains or patches $\Omega_j$ such that $\bigcup_{j=1}^{d} \Omega_j \supseteq \Omega$ with some mild overlap among the subdomains, we define the PUM by considering a partition of unity $\{w_j\}_{J=1}^{d}$ subordinated to the covering $\{\Omega_j\}_{j=1}^{d}$ such that $\sum_{j = 1}^d w_j( \boldsymbol{x})=1$, $\forall \boldsymbol{x} \in \Omega$. The weight $w_j: \Omega_j \rightarrow \RR$ is a compactly supported, nonnegative and continuous function with ${supp}(w_j)  \subseteq \Omega_j$. So the global RBF-PUM approximant assumes the form
	\begin{align} \label{intg}
	\tilde{u}( \boldsymbol{x})= \sum_{j=1}^{d} w_j ( \boldsymbol{x}) \tilde{u}_j( \boldsymbol{x} ), \quad \boldsymbol{x} \in \Omega.
	\end{align}
The local RBF approximant	$\tilde{u}_j:\Omega_j \rightarrow \RR$ is thus defined as
\begin{align} \label{loc_rbf}
	\tilde{u}_j( \boldsymbol{x}) = \sum_{i=1}^{N_j} c_i^j \phi_{\varepsilon} (||  \boldsymbol{x} - \boldsymbol{x}_i^j  ||_2),
\end{align}
where $N_j$ denotes the number of nodes in $\Omega_j$, i.e. the points $\boldsymbol{x}_i^j  \in X_{{N}_j} = X_N \cap \Omega_j$, $c_i^j$ is an unknown real coefficient, $||\cdot||_2$ indicates the Euclidean norm, and $\phi: \RR_{\geq 0}\to \RR$ is a RBF depending on a \textsl{shape parameter} $\varepsilon > 0$ such that $\phi_{\varepsilon}(||\boldsymbol{x}-\boldsymbol{z}||_2)=\phi(\varepsilon ||\boldsymbol{x}-\boldsymbol{z}||_2)$, $\forall \boldsymbol{x},\boldsymbol{z} \in \Omega$ (see \cite{fas07,liu18}). 

Thus, given the Laplace operator ${\cal L} = -\Delta$, we can define the Poisson problem with Dirichlet boundary conditions
\begin{align} \label{pde_poisson}
\begin{array}{rl}
-\Delta u(\boldsymbol{x}) = f(\boldsymbol{x}), & \quad  \boldsymbol{x} \in \Omega, \\
u(\boldsymbol{x}) = g(\boldsymbol{x}), & \quad  \boldsymbol{x} \in \partial\Omega.
\end{array}
\end{align}

The problem \eqref{pde_poisson} is then discretized on a global set of collocation points 
\begin{align*}
X_N = X_{N_i} \cup X_{N_b} = \{\boldsymbol{x}_1,\ldots,\boldsymbol{x}_N\} = \{\boldsymbol{x}_{i,1},\ldots,\boldsymbol{x}_{i,{N_i}}\} \cup \{x_{b,1},\ldots,x_{b,{N_b}}\},
\end{align*}
where $N_i$ and $N_b$ are the number of interior and boundary nodes, respectively. Assuming that the Poisson problem admits an approximate solution of the form \eqref{intg}, we get
\begin{align} \label{intg_poisson}
	\begin{array}{rl}
	\displaystyle{-\Delta\tilde{u}( \boldsymbol{x}_i)= -\sum_{j=1}^{d} \Delta \left(w_j ( \boldsymbol{x}_i) \tilde{u}_j( \boldsymbol{x}_i ) \right) = f(\boldsymbol{x}_i)}, &\quad \boldsymbol{x}_i \in \Omega\\
	\displaystyle{\tilde{u}( \boldsymbol{x}_i)= \sum_{j=1}^{d} w_j ( \boldsymbol{x}_i) \tilde{u}_j( \boldsymbol{x}_i ) = g(\boldsymbol{x}_i)}, &\quad \boldsymbol{x}_i \in \partial\Omega.
	\end{array}
\end{align}
The Laplace operator $\Delta$ can thus be expanded as
\begin{align} \label{lap_exp}
-\Delta \left(w_j(\boldsymbol{x}_i)\tilde{u}_j(\boldsymbol{x}_i)\right) = -\Delta w_j(\boldsymbol{x}_i)\tilde{u}_j(\boldsymbol{x}_i)-2\nabla w_j(\boldsymbol{x}_i) \cdot \nabla \tilde{u}_j(\boldsymbol{x}_i)-w_j(\boldsymbol{x}_i)\Delta\tilde{u}_j(\boldsymbol{x}_i), \quad \boldsymbol{x}_i \in \Omega.
\end{align}

Then we define the vector of local nodal values $\tilde{\boldsymbol{u}}_j = (\tilde{u}_j(\boldsymbol{x}_1^j), \ldots, \tilde{u}_j(\boldsymbol{x}_{N_j}^j))^T$, while the local coefficient vector $\boldsymbol{c}_j = (c_1^j,\ldots,c_{N_j}^j)^T$ is such that $\boldsymbol{c}_j = A_j^{-1} \tilde{\boldsymbol{u}}_j$. So we get
\begin{align} \label{diff_mat}
\Delta \tilde{\boldsymbol{u}}_j = A_j^{\Delta} A_j^{-1} \tilde{\boldsymbol{u}}_j, \qquad \nabla \tilde{\boldsymbol{u}}_j = A_j^{\nabla} A_j^{-1} \tilde{\boldsymbol{u}}_j,
\end{align}
where $A_j^{\Delta}$ and $A_j^{\nabla}$ are the matrices with entries
\begin{align*}
(A_j^{\Delta})_{ki} = \Delta \phi(||\boldsymbol{x}_k^j-\boldsymbol{x}_i^j||_2), \quad \mbox{and} \quad (A_j^{\nabla})_{ki} = \nabla \phi(||\boldsymbol{x}_k^j-\boldsymbol{x}_i^j||_2), \quad j = 1,\ldots, N_j.
\end{align*}
Moreover, we define the diagonal matrix 
\begin{align*}
W_j^{\Delta} = diag\left( \Delta w_j(\boldsymbol{x}_1^j),\ldots,\Delta w_j(\boldsymbol{x}_{N_j}^j)\right)
\end{align*}
associated with each subdomain, and similarly $W_j^{\nabla}$ and $W_j$. To obtain the discrete operator $L_j$, we have to differentiate \eqref{intg_poisson} by applying a product derivative rule and then use the relations in \eqref{diff_mat}. By using \eqref{lap_exp} and including the boundary conditions, we can express the discrete local Laplacian as
\begin{align*}
(L_j)_{ki}=
\left\{
\begin{array}{ll}
(\bar{L}_j)_{ki}, & \boldsymbol{x}_i^j \in \Omega,\\
\delta_{ki},      & \boldsymbol{x}_i^j \in \partial \Omega,
\end{array}
\right.
\end{align*} 
where $\delta_{ki}$ is the Kronecker delta and 
\begin{align*}
\bar{L}_j = \left( W_j^{\Delta} A_j + 2W_j^{\nabla}\cdot A_j^{\nabla}+W_j A_j^{\Delta}\right)A_j^{-1}.
\end{align*}
Finally, we obtain the global discrete operator by assembling the local matrices $L_j$ into the global matrix $L$ of entries
\begin{align*}
(L_j)_{ki} = \sum_{j=1}^{d} (L_j)_{\eta_{kj},\eta_{ij}}, \quad k,i=1,\ldots,N.
\end{align*} 
Thus, we have to solve the (global) linear system
\begin{align} \label{matcol}
L\boldsymbol{z} = \boldsymbol{u},
\end{align}
where $\boldsymbol{z}=(\tilde{u}(\boldsymbol{x}_1),\ldots,\tilde{u}(\boldsymbol{x}_N))^T$ and $\boldsymbol{u}=(u_1,\ldots, u_N)^T$ is defined through
\begin{align*}
u_i = \left\{
\begin{array}{ll}
f(\boldsymbol{x}_i), & \boldsymbol{x}_i \in \Omega, \\
g(\boldsymbol{x}_i), & \boldsymbol{x}_i \in \partial\Omega.
\end{array}
\right.
\end{align*}


\begin{remark}
In the PUM framework, an important computational issue consists in the organization of the collocation points among the subdomains $\Omega_j$, $j=1,\ldots,d$, since the computational efficiency is an essential requisite to fast construct the collocation matrix $L$ in \eqref{matcol} and accordingly find the solution of Poisson's equation via the RBF-PUM. The searching procedure used to localize and determine the points belonging to each subdomain $\Omega_j$ was first introduced to efficiently solve 2D and 3D interpolation problems \cite{cav16,cav18b}, and then further extended in higher dimensions \cite{all18}.  
\end{remark}


\section{Adaptive meshless refinement scheme} \label{amrs}

When seeking to construct an adaptive discretization method, one has to specify an error indicator and a refinement algorithm \cite{bab01}. To make an adaptive scheme effective, these two ingredients have to be appropriately selected so that each of them can take advantage of the potentiality of the other one, thus increasing the benefits when applying a numerical method to solve PDEs. However, as suggested in \cite{dri07}, a refinement consists generally in adding and/or removing collocation nodes but making use of the exact solution, an information that in applications is often not available. 


\subsection{Error indicators} 

In the work we thus aim to present a new adaptive scheme, which provides an error estimate without the exact solution is supposed to be known. This fact leads us to propose error indicators, which exclusively involve the solution of the collocation method or information to it related. In particular, we define two error indicators, referring explicitly to the iteration $k$ of the adaptive process. The first indicator is instead based on comparing the collocation RBF-PUM solution with the local RBF interpolant of the (approximate) collocation solution computed on $\Omega_j$ by computing the error on a set $Y^{(k)} = \{\boldsymbol{y}_1^{(k)}, \ldots, \boldsymbol{y}_{n^{(k)}}^{(k)}\}$ of test points. So the error indicator takes the form
\begin{align}
\label{indicator2}
e_i^{(k)} = \left|\tilde{u}(\boldsymbol{y}_i^{(k)}) - ({\cal I}\tilde{u})_{\Omega_j}(\boldsymbol{y}_i^{(k)})\right|, \qquad \boldsymbol{y}_i^{(k)} \in Y^{(k)},
\end{align}
where $\tilde{u}$ is the global RBF-PUM collocation solution computed on $X_{N^{(k)}}$ and $({\cal I}\tilde{u})_{\Omega_j}$ denotes the local RBF interpolant of $\tilde{u}$ constructed on the $N_j^{(k)}$ collocation points belonging to $\Omega_j$, $j=1,\ldots,d$ (cf. \eqref{loc_rbf}). The second indicator consists in comparing two collocation RBF-PUM solutions computed on two different data sets, a coarser set $X_{N_c^{(k)}}$ and a finer one $X_{N_f^{(k)}}$, such that $X_{N_c^{(k)}} \subseteq X_{N_f^{(k)}}$, $k=1,2,\ldots$, by evaluating then the error on the coarser set, i.e. 
\begin{align}
\label{indicator1}
e_i^{(k)} = \left|\tilde{u}_{X_{N_f^{(k)}}}(\boldsymbol{x}_i^{(k)}) - \tilde{u}_{X_{N_c^{(k)}}}(\boldsymbol{x}_i^{(k)})\right|, \qquad \boldsymbol{x}_i \in X_{N_c^{(k)}}.
\end{align}
Note that the PUM framework makes the indicator \eqref{indicator2} an ad-hoc strategy since we are clearly working on each PU subdomain. Therefore the latter is the one we will use in our extensive numerical tests. In fact, the indicator \eqref{indicator1} -- albeit its flexibility  -- is more expensive than \eqref{indicator2} from a computational standpoint. 


\subsection{Refinement algorithm} 

At first, we consider an initial set $X_{N^{(1)}}\equiv X_N = \{ \boldsymbol{x}_1^{(1)}, \ldots, \boldsymbol{x}_{N^{(1)}}^{(1)}\}$ of collocation points. This set is composed of two subsets which, in an adaptive algorithm, for $k=1,2,\ldots$, can iteratively be defined by a set $X_{N_i^{(k)}}=\{ \boldsymbol{x}_{i,1}^{(k)}, \ldots, \boldsymbol{x}_{i,N_i^{(k)}}^{(k)}\}$ of interior points and a set $X_{N_b^{(k)}}=\{ \boldsymbol{x}_{b,1}^{(k)}, \ldots, \boldsymbol{x}_{b,N_b^{(k)}}^{(k)}\}$ of boundary points, where $N_b^{(k)}= 4 \left\lceil \sqrt{N_i^{(k)}} + 2 \right\rceil - 4$.

Our adaptive meshless refinement scheme (AMRS) is based on the residual subsampling method proposed in \cite{dri07}. 
However, the mentioned paper is based on evaluating the residuals, i.e. supposing to know the exact solution of the Poisson problem. This fact highlights a significant weakness of the technique in \cite{dri07} because -- as specified earlier -- for instance in applications to science and engineering problems, it is essential for an adaptive discretization technique to combine a \emph{true} error indicator with a refinement procedure. Here, overall, we aims at proposing an adaptive scheme that allows one to be able to obtain reliable numerical results. The AMRS algorithm can thus be sketched as follows:
\begin{itemize}
	\item[\texttt{Step 1.}] We generate an initial discretization consisting of a set $X_{N^{(1)}}$ of grid (interior and boundary) collocation points and find the approximate solution on $X_{N^{(1)}}$ via RBF-PUM collocation scheme of the Poisson problem. 
	\item[\texttt{Step 2.}] We define a set $Y^{(1)}$ of Halton points as test points to compute/evaluate the error through the indicator \eqref{indicator2}. 
	\item[\texttt{Step 3.}] Fixed two tolerances $\tau_{\min} < \tau_{\max}$, the test points that exceed the threshold $\tau_{\max}$ are to be added among the collocation points, whereas the test point whose error is below the threshold $\tau_{\min}$ are removed along with their nearest point, i.e. we define the sets
	\begin{align*}
	Z_{T_{\max}^{(1)}} = \{\boldsymbol{y}_i^{(1)}\in Y^{(1)} \, : \, e_{i}^{(1)} > \tau_{\max}, \, i=1,\ldots,T_{\max}^{(1)}\} \quad \mbox{and} \quad
	Z_{T_{\min}^{(1)}} = \{\bar{\boldsymbol{x}}_i^{(1)}\in X_{N^{(1)}} \, : \, e_{i}^{(1)} < \tau_{\min}, \, i=1,\ldots,T_{\min}^{(1)}\},
	\end{align*}
	where $\bar{\boldsymbol{x}}_i^{(1)}$ is the nearest point to $\boldsymbol{y}_i^{(1)}$.
	\item[\texttt{Step 4.}] Iterating the process, for $k=2,3,\ldots$ we then construct a new set of discretization points
	\begin{align*}
	X_{N^{(k)}} = X_{N_i^{(k)}} \cup X_{N_b^{(k)}}, \quad	\mbox{with} \quad	X_{N_i^{(k)}} = (X_{N_i^{(k-1)}} \cup Z_{T_{\max}^{(k-1)}}) \backslash Z_{T_{\min}^{(k-1)}},
	\end{align*}
	and repeat the procedure as in \texttt{Step 1-3}, also defining a new set $Y^{(k)}$ of test points.
	\item[\texttt{Step 5.}] We stop when the set $Z_{T_{\min}^{(k)}}$ is empty.
\end{itemize}
This adaptive refinement scheme follows therefore the common paradigm of solve-estimate-refine/coarsen till a stopping criterion is satisfied.


\section{Numerical experiments} \label{num_exp}

In this section we illustrate the performance of our adaptive RBF-PUM algorithm, which is implemented in \textsc{Matlab}. All the results are carried out on a laptop with an Intel(R) Core(TM) i7-4500U CPU 1.80 GHz processor and 4GB RAM. In the following we focus on a wide series of experiments, which concern the numerical solution of Poisson problems via RBF-PUM collocation as described in Ssection \ref{RBF-PUM}. In particular, we test our AMRS technique, which involves the use of error indicator \eqref{indicator2} with the related adaptive refinement, as outlined in Section \ref{amrs}.

\begin{figure}[ht!]
\centering
\parbox{7.5cm}{\centering
\includegraphics[scale=0.4]{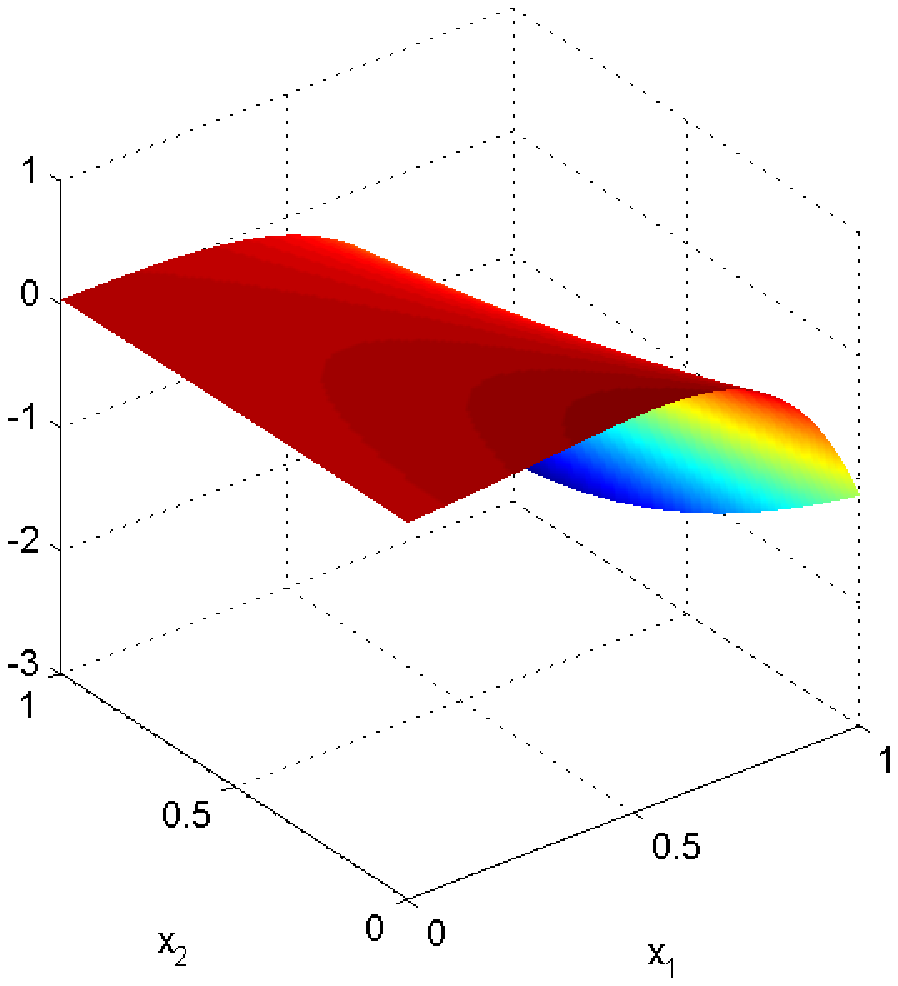}
\captionof*{figure}{$u_1(x_1,x_2)$}}
\quad
\parbox{7.5cm}{\centering
\includegraphics[scale=0.4]{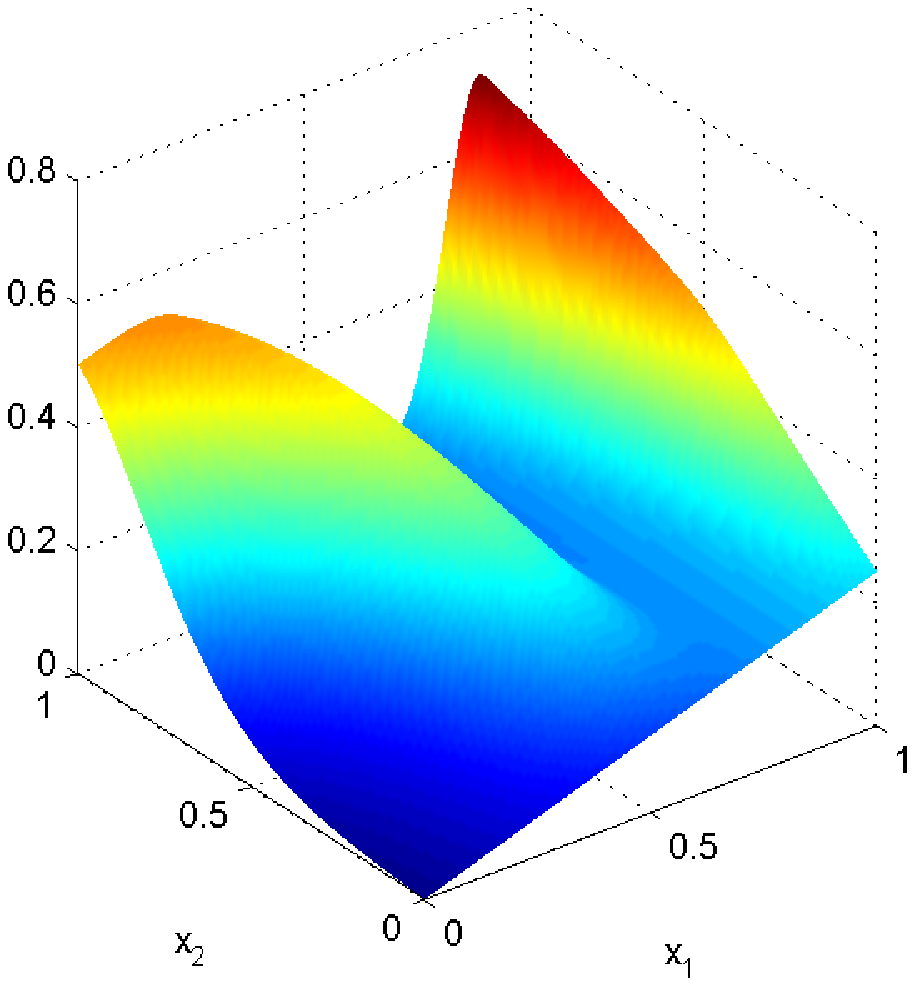}
\captionof*{figure}{$u_2(x_1,x_2)$}}\\
\parbox{7.5cm}{\centering
\includegraphics[scale=0.4]{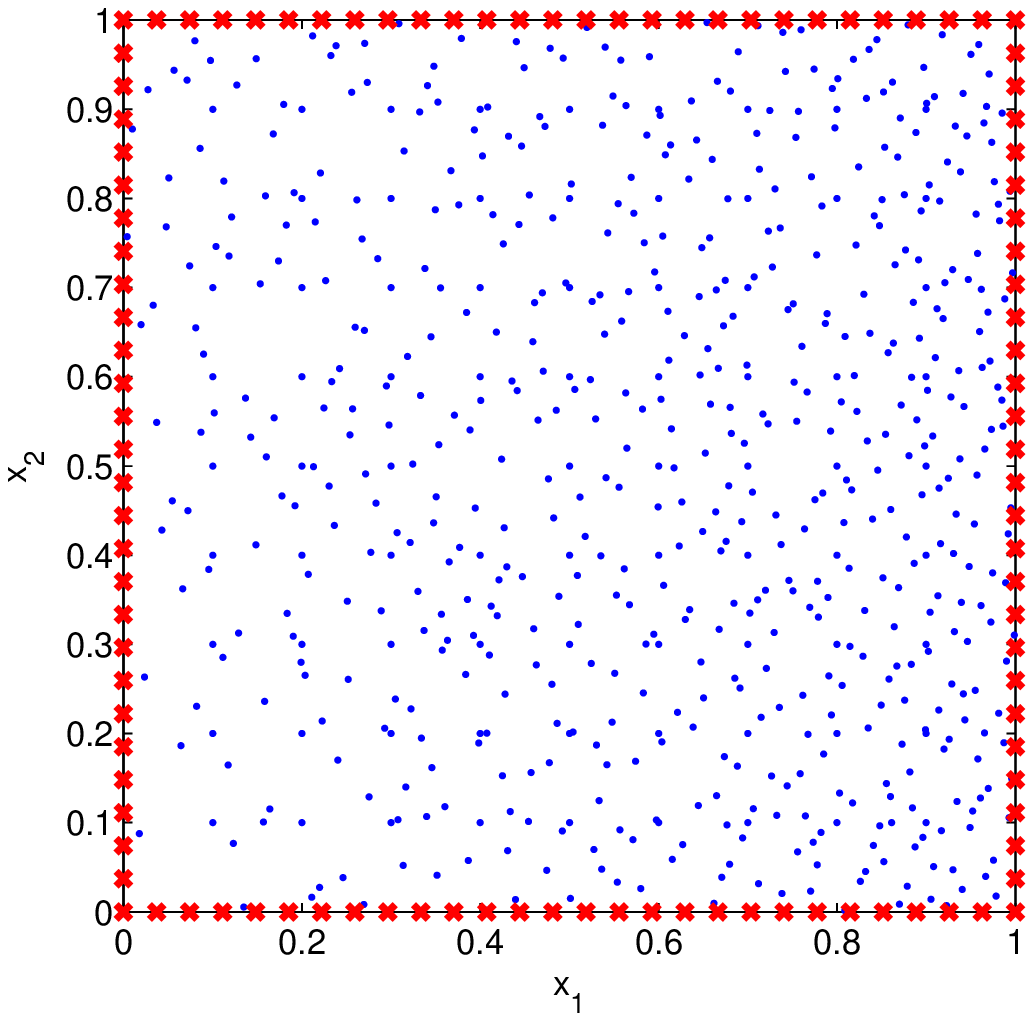}
\captionof*{figure}{$u_1$: final data}}
\quad
\parbox{7.5cm}{\centering
\includegraphics[scale=0.4]{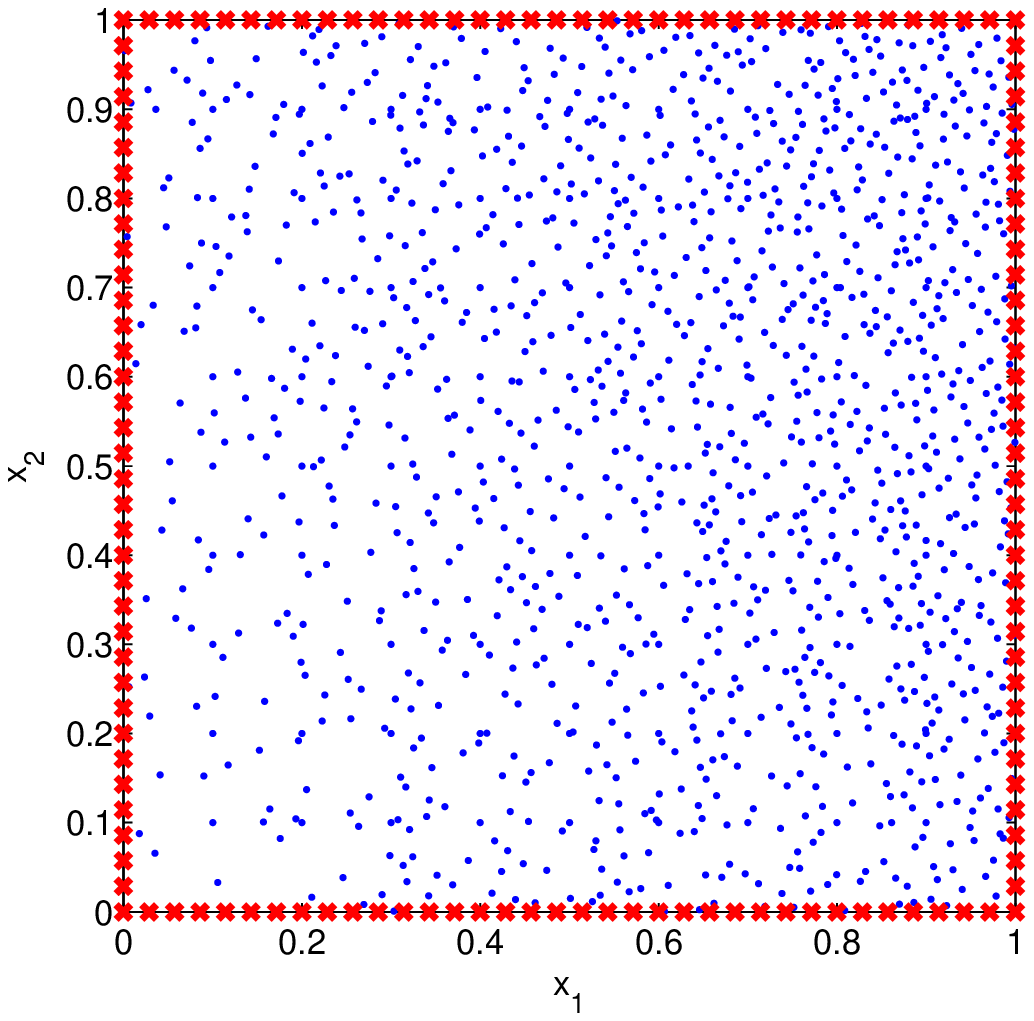}
\captionof*{figure}{$u_2$: final data}}
\caption{Exact solutions of Poisson problems (top) and final data distribution in case (a) obtained via AMRS for RBF-PUM collocation using M6 for $\varepsilon = 3$ (bottom).}
\label{fig:poisson_sol}
\end{figure}

In these tests we thus show the results obtained by applying the RBF-PUM scheme using the Mat$\acute{\text{e}}$rn function M6, i.e. $\phi_{\varepsilon}(r)=\exp({-\varepsilon r}) (\varepsilon^3 r^3 + 6\varepsilon^2r^2+15\varepsilon r+15)$, as local approximant in \eqref{loc_rbf} with shape parameter $\varepsilon = 3$ and the compactly supported Wendland function W2 as localizing function of Shepard's weight. We analyze the behavior of the adaptive algorithm by considering a number of test problems, where an adaptive refinement is known to be of great advantage. So, referring to \eqref{pde_poisson}, we focus on two Poisson problems defined on $\Omega=[0,1]\times[0,1]$, whose exact solutions are given by
\begin{align*}
u_1(x_1,x_2) = \displaystyle{\frac{1}{20} \exp(4x_1)\cos(2x_1+x_2)}, \qquad
u_2(x_1,x_2) = \displaystyle{\frac{1}{2}x_2 \left[\cos(4x_1^2+x_2^2-1)\right]^4+\frac{1}{4}x_1}.
\end{align*}
In Figure \ref{fig:poisson_sol} (top) we show a graphical representation of such analytic solutions. However, we remark that several numerical experiments (not reported here for shortness) have been carried out using other test problems and the results show a uniform behavior.

As a measure of the quality/accuracy of the results, we compute the Maximum Absolute Error (MAE) and the Root Mean Square Error (RMSE) on a grid of $N_{eval} = 40 \times 40$ evaluation points, i.e.
\begin{align} \label{mae_rmse_errors}
	\mbox{MAE} = \max_{1 \leq i \leq N_{eval}} |u(\boldsymbol{z}_i) - \tilde{u}(\boldsymbol{z}_i)|, \qquad	\mbox{RMSE} = \sqrt{\frac{1}{N_{eval}}\sum_{i=1}^{N_{eval}} |u(\boldsymbol{z}_i) - \tilde{u}(\boldsymbol{z}_i)|^2}.
\end{align}
Moreover, in order to analyze the stability of the method, we evaluate the Condition Number (CN) of the sparse collocation matrix $L$ in \eqref{matcol} by using the \textsc{Matlab} command \texttt{condest}. As regards, instead, the efficiency of our algorithm, we report the CPU times computed in seconds.

The main target of our numerical analysis is thus addressed to highlight the performance of the adaptive scheme. In particular, in Tables \ref{tab:1}--\ref{tab:1b} we report the results obtained by assuming $(\tau_{\min},\tau_{\max})=(10^{-8},10^{-5})$ and starting from $N=121$ collocation points, which consist of $N_b=40$ grid points and: (a) $N_i = 81$ grid points; (b) $N_i = 81$ Halton points. Moreover, we also indicate the final number $N_{tot}$ of collocation points required to reach the min/max thresholds. For brevity, in Figure \ref{fig:poisson_sol} (bottom) we plot the final data sets only focusing on the case (a).

Analyzing the numerical results, we can observe that the adaptive scheme works well increasing the number of points in the regions where the surface of the solution is characterized by significant variations. Then, focusing on the stability, we report the CN that in all problems faced has a order of magnitude between $10^{+6}$ and $10^{+9}$. Note that the latter is quite lower than that observed in traditional RBF-based collocation methods, where the CN is often larger than $10^{+18}$ \cite{fas07}. Finally, as to efficiency of the iterative algorithm, execution times shown in Tables \ref{tab:1}--\ref{tab:1b} highlight as the numerical scheme stops (after reaching the given thresholds) in few seconds.

\begin{table}[ht!]
{\small
\begin{center} 
{
\begin{tabular}{cccccc} 
\toprule
Poisson & $N_{tot}$ & MAE & RMSE & CN & time   \\ \midrule
$u_1$  & $755$ & $1.00$e${-4}$ & $3.46$e${-5}$ & $6.58$e${+06}$ & $4.6$  \\ \midrule
$u_2$  & $1411$ & $1.58$e${-4}$ & $1.71$e${-5}$ & $4.04$e${+08}$ & $8.6$  \\ \midrule
\end{tabular}
}
\caption{Results in (a) obtained via AMRS for RBF-PUM collocation using M6 for $\varepsilon = 3$.} 
\label{tab:1}
\end{center}
}
\end{table}

\begin{table}[ht!]
{\small
\begin{center} 
{
\begin{tabular}{cccccc} \toprule
Poisson & $N_{tot}$ & MAE & RMSE & CN & time   \\ \midrule
$u_1$  & $746$ & $1.21$e${-4}$ & $3.45$e${-5}$ & $9.56$e${+06}$ & $4.0$  \\ \midrule
$u_2$  & $1452$ & $1.13$e${-4}$ & $1.72$e${-5}$ & $1.31$e${+08}$ & $8.8$  \\ \midrule
\end{tabular}
}
\caption{Results in (b) obtained via AMRS for RBF-PUM collocation using M6 for $\varepsilon = 3$.} 
\label{tab:1b}
\end{center}
}
\end{table}

\section*{Acknowledgments}
The authors acknowledge support from the Department of Mathematics \lq\lq Giuseppe Peano\rq\rq\ of the University of Torino via 2016-2017 projects \lq\lq Numerical and computational methods for applied sciences\rq\rq\ and \lq\lq Multivariate approximation and efficient algorithms with applications to algebraic, differential and integral problems\rq\rq. Moreover, this work was partially supported by GNCS--INdAM 2018 project. This research has been accomplished within RITA (Research ITalian network on Approximation).





\end{document}